\documentclass[11pt]{article}
\usepackage{amsfonts}
\usepackage{amsfonts}
\usepackage{amsfonts}
\usepackage{mathrsfs}
\usepackage{bm}
\usepackage{cite}
\usepackage[colorlinks,linkcolor=blue,citecolor=blue]{hyperref}
\usepackage{amssymb,amsmath} % font used for R in Real numbers
 \allowdisplaybreaks

\catcode`!=11
\let\!int\int \def\int{\displaystyle\!int}
\let\!lim\lim \def\lim{\displaystyle\!lim}
\let\!sum\sum \def\sum{\displaystyle\!sum}
\let\!sup\sup \def\sup{\displaystyle\!sup}
\let\!inf\inf \def\inf{\displaystyle\!inf}
\let\!cap\cap \def\cap{\displaystyle\!cap}
\let\!max\max \def\max{\displaystyle\!max}
\let\!min\min \def\min{\displaystyle\!min}
\let\!frac\frac \def\frac{\displaystyle\!frac}
\catcode`!=12

\let\oldsection\section
\renewcommand\section{\setcounter{equation}{0}\oldsection}

\allowdisplaybreaks
\def\pf{\it{Proof.}\rm\quad}

\def\N{\mathbb{N}}\def\Z{\mathbb{Z}}

\newtheorem{thm}{Theorem}[section]
\newtheorem{lem}[thm]{Lemma}

\newtheorem{re}{Remark}[section]
\begin{document}
%%%%%%%%%%%%%%%%%%%% title %%%%%%%%%%%%%%%%%%%%%%%%%%%%%%%%%%%%%%%%%%%%%%%%
\title {\bf  Computation and theory of Euler sums of generalized hyperharmonic numbers}
\author{
{Ce Xu\thanks{Corresponding author. Email: 15959259051@163.com (C. Xu)}}\\[1mm]
\small  School of Mathematical Sciences, Xiamen University\\
\small Xiamen
361005, P.R. China}
\date{}
\maketitle \noindent{\bf Abstract } Recently, Dil and Boyadzhiev \cite{AD2015} proved an explicit formula for the sum of multiple harmonic numbers whose indices are the sequence $\left( {{{\left\{ 0 \right\}}_r},1} \right)$. In this paper we show that the sums of multiple harmonic numbers whose indices are the sequence $\left( {{{\left\{ 0 \right\}}_r,1};{{\left\{ 1 \right\}}_{k-1}}} \right)$ can be expressed in terms of (multiple) zeta values, multiple harmonic numbers and Stirling numbers of the first kind, and give an explicit formula.
\\[2mm]
\\[2mm]
\tableofcontents
\section{Introduction}
Let $s_1,\ldots,s_m$ be positive integers. The classical multiple harmonic numbers (MHNs) and multiple harmonic star numbers (MHSNs) are defined by the partial sums (see \cite{KP2013,Xu2017})
\[{\zeta _n}\left( {{s_1},{s_2}, \cdots ,{s_m}} \right): = \sum\limits_{n \ge {n_1}  >  \cdots  > {n_m} \ge 1} {\frac{1}{{n_1^{{s_1}}n_2^{{s_2}} \cdots n_m^{{s_m}}}}} ,\tag{1.1}\]
\[{\zeta_n ^ \star }\left( {{s_1},{s_2}, \cdots ,{s_m}} \right): = \sum\limits_{n \ge {n_1}  \ge  \cdots  \ge {n_m} \ge 1} {\frac{1}{{n_1^{{s_1}}n_2^{{s_2}} \cdots n_m^{{s_m}}}}},\tag{1.2}\]
when $n<m$, then ${\zeta _n}\left( {{s_1},{s_2}, \cdots ,{s_m}} \right)=0$, and ${\zeta _n}\left(\emptyset \right)={\zeta^\star _n}\left(\emptyset \right)=1$. The limit cases of MHNs and MHSNs give rise to multiple zeta values (MZVs) and multiple zeta star values (MZSVs) (see \cite{KO2010,OZ2001,KP2013,Xu2017,DZ2012}):
\begin{align*}
&\zeta \left( {{s_1},{s_2}, \cdots ,{s_m}} \right) = \mathop {\lim }\limits_{n \to \infty } \zeta \left( {{s_1},{s_2}, \cdots ,{s_m}} \right),\\
&\zeta _{}^ \star \left( {{s_1},{s_2}, \cdots ,{s_m}} \right) = \mathop {\lim }\limits_{n \to \infty } \zeta _{}^ \star \left( {{s_1},{s_2}, \cdots ,{s_m}} \right)
\end{align*}
defined for $s_2,\ldots,s_m\geq 1$ and $s_1\geq 2$ to ensure convergence of the series.
For non-negative integers $s_1,\ldots,s_{m+k}$, we define the following a generalized multiple harmonic numbers
\[{\zeta _n}\left( {{s_1}, \cdots ,{s_m};{s_{m + 1}}, \cdots ,{s_{m + k}}} \right): = \sum\limits_{\scriptstyle 0 < {n_{m + k}} <  \cdots  < {n_{m+1}} \hfill \atop
  \scriptstyle  <n_m\le {n_{m - 1}} \cdots  \le {n_1} \le n \hfill} {\frac{1}{{n_1^{{s_1}} \cdots n_{m + k}^{{s_{m + k}}}}}}.\tag{1.3} \]
Obviously, if $m=0$ or $k=0$ in (1.3) and $s_j\in \N:=\{1,2,\cdots\}$, then
\begin{align*}
&{\zeta _n}\left( {\emptyset ;{s_1}, \cdots ,{s_k}} \right) = {\zeta _n}\left( {{s_1}, \cdots ,{s_k}} \right),\\
&{\zeta _n}\left( {{s_1},\cdots ,{s_m};\emptyset } \right) = \zeta _n^ \star \left( {{s_1}, \cdots ,{s_m}} \right).
\end{align*}
There are a lot of recent contributions on MZVs and MZSVs (for example, see \cite{KO2010,OZ2001,KP2013,Xu2017,DZ2012}).
The earliest results on MZVs or MZSVs are due to Euler who elaborated a method to reduce double sums $\zeta \left( {{s_1},{s_2}} \right)$ (also called linear Euler sums \cite{FS1998,X2016}) of small weight to certain rational linear combinations of products of zeta values. In \cite{FS1998}, Flajolet and Salvy introduced the following generalized series
\[{S_{{\bf S},q}} := \sum\limits_{n = 1}^\infty  {\frac{{H_n^{\left( {{s_1}} \right)}H_n^{\left( {{s_2}} \right)} \cdots H_n^{\left( {{s_r}} \right)}}}
{{{n^q}}}},\tag{1.4}\]
which is called the generalized (nonlinear) Euler sums. Here ${\bf S}:=(s_1,s_2,\ldots,s_r)\ (r,s_i\in \N, i=1,2,\ldots,r)$ with $s_1\leq s_2\leq \ldots\leq s_r$ and $q\geq 2$. The quantity $w:={s _1} +  \cdots  + {s _r} + q$ is called the weight and the quantity $r$ is called the degree. The notation $H_n^{\left( p \right)}$ denotes the ordinary harmonic numbers defined by
\[H_n^{\left( p \right)}\equiv \zeta_n(p): = \sum\limits_{j = 1}^n {\frac{1}{{{j^p}}}},\ p,n\in \N. \]
It has been discovered in the course of the years that many nonlinear Euler sums admit expressions involving finitely "zeta values", that is say values of the Riemann zeta function,
$$\zeta(s):=\sum\limits_{n = 1}^\infty {\frac {1}{n^{s}}},\Re(s)>1.$$
with positive integer arguments, and linear Euler sums.
The relationship between the values of the Riemann zeta values and Euler sums has been studied by many authors. For details and historical introductions, please see \cite{BBG1994,BBG1995,BBGP1996,BZB2008,BG1996,AD2009,FS1998,M2014,S2015,S2012,Xu2016,XZZ2016,X2016} and references therein.

From \cite{CG1996,AD2008,AD2015,M2010,GS2017}, we know that the classical hyperharmonic numbers are defined by
\[h_n^{\left( m \right)}: = \sum\limits_{1 \le {n_m} \le  \cdots  \le {n_1} \le n} {\frac{1}{{{n_m}}}}  = \zeta _n^ \star \left( {{{\left\{ 0 \right\}}_{m - 1}},1} \right).\tag{1.5}\]
In \cite{X2017}, we define the generalized hyperharmonic numbers $h_n^{\left( m \right)}(k)$ by
\[h_n^{\left( m \right)}\left( k \right): = \sum\limits_{\scriptstyle 1 \le {n_{m + k - 1}} <  \cdots  < {n_{m+1}} \hfill \atop
  \scriptstyle  <n_m\le {n_{m - 1}} \le  \cdots  \le {n_1} \le n \hfill} {\frac{1}{{{n_m}{n_{m + 1}} \cdots {n_{m + k - 1}}}}}  = {\zeta _n}\left( {{{\left\{ 0 \right\}}_{m - 1}},1;{{\left\{ 1 \right\}}_{k - 1}}} \right),\tag{1.6}\]
where $m,k\in \N$ \ (The notation $\{ \}_p$ means that the sequence in the bracket is repeated $p$ times). In this paper, we prove the result: for positive integers $m$ and $k$, the Euler-type sums with hyperharmonic numbers
\[S\left( {k,m;p} \right): = \sum\limits_{n = 1}^\infty  {\frac{{h_n^{\left( m \right)}\left( k \right)}}{{{n^p}}}} \;\;\left(p\geq m+1 \right)\]
are related to the multiple zeta values, multiple harmonic numbers and Stirling numbers of the first
kind. For $k=1,2,3$, the above results have been proved in Dil et al.\cite{AD2015} and our paper \cite{X2017}.
The purpose of the present paper is to prove the following theorems.
\begin{thm} For integers $k\in \mathbb{N}$ and $r\in \mathbb{N}_0$ with $r+2\leq p\in \mathbb{N}$, then the following identity holds:
\[S\left( {k,r + 1;p} \right) = \frac{1}{{r!}}\sum\limits_{l = 1}^{r + 1} {\left[ {\begin{array}{*{20}{c}}
   {r + 1}  \\
   l  \\
\end{array}} \right]} \sum\limits_{\scriptstyle i + j = k \hfill \atop
  \scriptstyle i,j \ge 0 \hfill} {{{\left( { - 1} \right)}^i}{\zeta_r ^ \star }\left( {{{\left\{ 1 \right\}}_i}} \right){U_{j,r}}\left( {p + 1 - l} \right)}, \tag{1.7}\]
where $\left[ {\begin{array}{*{20}{c}}
   n  \\
   k  \\
\end{array}} \right]$  denotes the (unsigned) Stirling number of the first kind, which is defined by \cite{CG1996,L1974}
\[n!x\left( {1 + x} \right)\left( {1 + \frac{x}{2}} \right) \cdots \left( {1 + \frac{x}{n}} \right) = \sum\limits_{k = 0}^n {\left[ {\begin{array}{*{20}{c}}
   {n + 1}  \\
   {k + 1}  \\
\end{array}} \right]{x^{k + 1}}}, \tag{1.8} \]
with $\left[ {\begin{array}{*{20}{c}}
   n  \\
   k  \\
\end{array}} \right]=0$, if $n<k$ and $\left[ {\begin{array}{*{20}{c}}
   n  \\
   0  \\
\end{array}} \right]=\left[ {\begin{array}{*{20}{c}}
   0  \\
   k  \\
\end{array}} \right]=0,\ \left[ {\begin{array}{*{20}{c}}
   0  \\
   0  \\
\end{array}} \right]=1$, or equivalently, by the generating function:
 \[{\log ^k}\left( {1 - x} \right) = {\left( { - 1} \right)^k}k!\sum\limits_{n = 1}^\infty  {\left[ {\begin{array}{*{20}{c}}
   n  \\
   k  \\
\end{array}} \right]\frac{{{x^n}}}{{n!}}} ,\;x \in \left[ { - 1,1} \right).\tag{1.9} \]
and
\[{U_{j,r}}\left( p \right): = \sum\limits_{n = 1}^\infty  {\frac{{{\zeta _{n + r}}\left( {{{\left\{ 1 \right\}}_j}} \right)}}{{{n^p}}}} .\tag{1.10}\]
\end{thm}
\begin{thm}\label{thm2} For integers $m>0,p>1$ and $r\in \N_0:=\N\cup \{0\}$, we have
\begin{align*}
&{U_{m,r}}\left( p \right) = \zeta \left( {p,{{\left\{ 1 \right\}}_m}} \right) + \zeta \left( {p + 1,{{\left\{ 1 \right\}}_{m - 1}}} \right)\\
& + {\left( { - 1} \right)^{p - 1}}\sum\limits_{1 \le {i_m} <  \cdots  < {i_1} \le r} {\sum\limits_{l = 1}^m {\prod\limits_{\scriptstyle a = 1 \hfill \atop
  \scriptstyle a \ne l \hfill}^m {\frac{{{H_{{i_l}}}}}{{\left( {{i_a} - {i_l}} \right)i_l^p}}} } } \\& + \sum\limits_{1 \le {i_m} <  \cdots  < {i_1} \le r} {\sum\limits_{b = 1}^{p - 1} {\sum\limits_{l = 1}^m {\prod\limits_{\scriptstyle a = 1 \hfill \atop
  \scriptstyle a \ne l \hfill}^m {\frac{{\zeta \left( {p + 1 - b} \right)}}{{\left( {{i_a} - {i_l}} \right)i_l^b}}{{\left( { - 1} \right)}^{b - 1}}} } } } \\
  &+ {\left( { - 1} \right)^{p-1}}\sum\limits_{j = 1}^{m - 1} {\sum\limits_{1 \le {i_j} <  \cdots  < {i_1} \le r} {\sum\limits_{l = 1}^j {\prod\limits_{\scriptstyle a = 1 \hfill \atop
  \scriptstyle a \ne l \hfill}^j {\frac{{\zeta \left( {m - j + 1} \right) + \zeta _{{i_l}}^ \star \left( {{{\left\{ 1 \right\}}_{m - j + 1}}} \right) - \frac{{\zeta _{{i_l}}^ \star \left( {{{\left\{ 1 \right\}}_{m - j}}} \right)}}{{{i_l}}}}}{{\left( {{i_a} - {i_l}} \right)i_l^p}}} } } } \\
&+ \sum\limits_{j = 1}^{m - 1} {\sum\limits_{1 \le {i_j} <  \cdots  < {i_1} \le r} {\sum\limits_{b = 1}^{p - 1} {\sum\limits_{l = 1}^j {\prod\limits_{\scriptstyle a = 1 \hfill \atop
  \scriptstyle a \ne l \hfill}^j {\frac{{\zeta \left( {p + 1 - b,{{\left\{ 1 \right\}}_{m - j}}} \right) + \zeta \left( {p + 2 - b,{{\left\{ 1 \right\}}_{m - j - 1}}} \right)}}{{\left( {{i_a} - {i_l}} \right)i_l^b}}{{\left( { - 1} \right)}^{b - 1}}} } } } } ,\tag{1.11}
\end{align*}
where $H_{i_l}$ is harmonic number.
\end{thm}

It is clear that the Theorem \ref{thm2} implies that the sums ${U_{m,r}}\left( p \right)$ can be expressed in terms of series of Riemann zeta function and harmonic numbers.
\section{Some Lemmas and Theorems}
To prove the Theorem 1.1 and Theorem 1.2, we need the following lemmas.
\begin{lem}(\cite{X2017})
For positive integers $n$ and $k$, then the following identity holds:
\[\left[ {\begin{array}{*{20}{c}}
   n  \\
   k  \\
\end{array}} \right] = \left( {n - 1} \right)!{\zeta _{n - 1}}\left( {{{\left\{ 1 \right\}}_{k - 1}}} \right).\tag{2.1}\]
\end{lem}
\begin{lem} (\cite{X2017}) For positive integers $m,n$ and $k$, we have the recurrence relation
\[h_n^{\left( m \right)}\left( k \right) = \frac{{{{\left( { - 1} \right)}^{k - 1}}}}{k}\sum\limits_{i = 0}^{k - 1} {{{\left( { - 1} \right)}^i}h_n^{\left( m \right)}\left( i \right)\left\{ {H_{m + n - 1}^{\left( {k - i} \right)} - H_{m - 1}^{\left( {k - i} \right)}} \right\}} ,\tag{2.2}\]
where \[h_n^{\left( m \right)}\left( 0 \right) = \left( {\begin{array}{*{20}{c}}
   {m + n - 1}  \\
   {m - 1}  \\
\end{array}} \right).\]
\end{lem}
\begin{lem} (\cite{Xu2017}) For positive integers $m$ and $n$, then the recurrence relation holds:
\[{{\bar B}_m}\left( n \right) = \frac{{{{\left( { - 1} \right)}^{m - 1}}}}{m}\sum\limits_{i = 0}^{m - 1} {{{\left( { - 1} \right)}^i}{{\bar B}_i}\left( n \right){X_n}\left( {m - i} \right)},\tag{2.3}\]
where
\[{X_n}\left( m \right) := \sum\limits_{i = 1}^n {x_i^m},\ x_i\in \mathbb{C},\]
\[{\bar B}_m(n) := \underset{k_1 = 1}{\overset{n}{\sum}}x_{k_1}\underset{k_2 = 1}{\overset{k_1-1}{\sum}}x_{k_2}\cdots\underset{k_m = 1}{\overset{k_{m-1}-1}{\sum}}x_{k_m},{\bar B}_0(n) = 1.\]
\end{lem}
\begin{lem} (\cite{XZZ2016,X2017})
For integers $k\in \N$ and $p\in \mathbb{N} \setminus \{1\}:=\{2,3,\ldots\}$,  then the following identity holds:
\[\left( {p - 1} \right)!\sum\limits_{n = 1}^\infty  {\frac{{\left[ {\begin{array}{*{20}{c}}
   {n + 1}  \\
   p  \\
\end{array}} \right]}}{{n!n\left( {n + k} \right)}}}  = \frac{1}{k}\left\{ {\left( {p - 1} \right)!\zeta (p) + \frac{{{Y_p}\left( k \right)}}{p} - \frac{{{Y_{p - 1}}\left( k \right)}}{k}} \right\},\tag{2.4}\]
where ${Y_k}\left( n \right) = {Y_k}\left( {{H _n},1!{H^{(2)} _n},2!{H^{(3)}_n}, \cdots ,\left( {r - 1} \right)!{H^{(r)} _n}, \cdots } \right)$, ${Y_k}\left( {{x_1},{x_2}, \cdots } \right)$ stands for the complete exponential Bell polynomial defined by (see \cite{L1974})
\[\exp \left( {\sum\limits_{m \ge 1}^{} {{x_m}\frac{{{t^m}}}{{m!}}} } \right) = 1 + \sum\limits_{k \ge 1}^{} {{Y_k}\left( {{x_1},{x_2}, \cdots } \right)\frac{{{t^k}}}{{k!}}}.\tag{2.5}\]
\end{lem}
Noting that, in \cite{X2016}, we find the relation
\[\zeta _n^ \star \left( {{{\{ 1\} }_m}} \right) = \frac{1}{{m!}}{Y_m}\left( n \right),n,m\in \N_0.\tag{2.6}\]
\begin{lem} For positive integers $m$ and $n$, then
\[\sum\limits_{\scriptstyle i + j = m \hfill \atop
  \scriptstyle i,j \ge 0 \hfill} {{{\left( { - 1} \right)}^i}{{\bar A}_i}\left( n \right){A_j}\left( n \right)}  = 0,\tag{2.7}\]
  where
\[{{\bar A}_m}\left( n \right): = \sum\limits_{1 \le {k_m} <  \cdots  < {k_1} \le n} {{a_{{k_1}}} \cdots {a_{{k_m}}}} ,\;{a_k} \in \mathbb{C},\]
\[{A_m}\left( n \right) := \sum\limits_{1 \le {k_m} \le  \cdots  \le {k_1} \le n} {{a_{{k_1}}} \cdots {a_{{k_m}}}} ,\;{a_k} \in \mathbb{C}.\]
For convenience, we set ${{\bar A}_0}\left( n \right) = {A_0}\left( n \right) = 1$. If $n<m$, we let ${{\bar A}_m}\left( n \right) = 0$.
\end{lem}
\pf By a direct calculation, the following identities are easily derived
\begin{align*}
&\prod\limits_{i = 1}^n {{{\left( {1 - {a_i}t} \right)}^{ - 1}}}  = \sum\limits_{m = 0}^\infty  {{A_m}\left( n \right){t^m}} ,\\
&\prod\limits_{i = 1}^n {\left( {1 + {a_i}t} \right)}  = \sum\limits_{m = 0}^\infty  {{{\bar A}_m}\left( n \right){t^m}} .
\end{align*}
Hence, by using Cauchy product of power series, we have
\begin{align*}
1 &= \prod\limits_{i = 1}^n {{{\left( {1 - {a_i}t} \right)}^{ - 1}}} \prod\limits_{i = 1}^n {\left( {1 - {a_i}t} \right)} \\
& = \left( {\sum\limits_{m = 0}^\infty  {{A_m}\left( n \right){t^m}} } \right)\left( {\sum\limits_{m = 0}^\infty  {{{\left( { - 1} \right)}^m}{{\bar A}_m}\left( n \right){t^m}} } \right)\\
& = \sum\limits_{m = 0}^\infty  {\left\{ {\sum\limits_{\scriptstyle i + j = m \hfill \atop
  \scriptstyle i,j \ge 0 \hfill} {{{\left( { - 1} \right)}^i}{{\bar A}_i}\left( n \right){A_j}\left( n \right)} } \right\}{t^m}} .
\end{align*}
Thus, comparing the coefficients of $t^m$ in above equation, we obtain the formula (2.7). The proof of Theorem 2.3 is finished. \hfill$\square$\\
The above lemmas will be useful in the development of the main theorems. Next, we give some important theorems and it's proofs by using these lemmas.
\begin{thm} For integers $r\geq 0$ and $m,n>1$, then
\[\sum\limits_{1 \le {k_m} <  \cdots  < {k_1} \le n} {\frac{1}{{\left( {{k_1} + r} \right) \cdots \left( {{k_m} + r} \right)}}}  = \sum\limits_{\scriptstyle i + j = m \hfill \atop
  \scriptstyle i,j \ge 0 \hfill} {{{\left( { - 1} \right)}^i}\zeta _r^ \star \left( {{{\left\{ 1 \right\}}_i}} \right){\zeta _{n + r}}\left( {{{\left\{ 1 \right\}}_j}} \right)} ,\tag{2.8}\]
  where \[\zeta _0^ \star \left( \emptyset  \right) = {\zeta _0}\left( \emptyset  \right) = 1,\;\zeta _0^ \star \left( {{{\left\{ 1 \right\}}_i}} \right) = {\zeta _0}\left( {{{\left\{ 1 \right\}}_j}} \right) = 0,\;\left( {i,j \ge 1} \right).\]

\end{thm}
\pf The proof is by induction on $m$. For $m=1$ we have $\sum\limits_{1 \le {k_1} \le n} {\frac{1}{{{k_1} + r}}}  = {\zeta _{n + r}}\left( 1 \right) - \zeta _r^ \star \left( 1 \right)$, and the formula is true. For $m>1$ we proceed as follow. Let
\[{\zeta _n}\left( {{s_1},{s_2}, \cdots ,{s_m}\left| {r + 1} \right.} \right): = \sum\limits_{1 \le {k_m} <  \cdots  < {k_1} \le n} {\frac{1}{{{{\left( {{k_1} + r} \right)}^{{s_1}}} \cdots {{\left( {{k_m} + r} \right)}^{{s_m}}}}}} ,\tag{2.9}\]
\[{\zeta _n}\left( {\emptyset \left| {r + 1} \right.} \right) = 1.\]
Then by the definition (2.9) and the induction hypothesis, we have that
\begin{align*}
{\zeta _n}\left( {{{\left\{ 1 \right\}}_{m + 1}}\left| {r + 1} \right.} \right) &= \sum\limits_{k = 1}^n {\frac{{{\zeta _{k - 1}}\left( {{{\left\{ 1 \right\}}_m}\left| {r + 1} \right.} \right)}}{{k + r}}} \\
&= \sum\limits_{k = 1}^n {\frac{1}{{k + r}}\sum\limits_{\scriptstyle i + j = m \hfill \atop
  \scriptstyle i,j \ge 0 \hfill} {{{\left( { - 1} \right)}^i}\zeta _r^ \star \left( {{{\left\{ 1 \right\}}_i}} \right){\zeta _{k + r - 1}}\left( {{{\left\{ 1 \right\}}_j}} \right)} } \\
& = \sum\limits_{\scriptstyle i + j = m \hfill \atop
  \scriptstyle i,j \ge 0 \hfill} {{{\left( { - 1} \right)}^i}\zeta _r^ \star \left( {{{\left\{ 1 \right\}}_i}} \right)\sum\limits_{k = 1}^n {\frac{{{\zeta _{k + r - 1}}\left( {{{\left\{ 1 \right\}}_j}} \right)}}{{k + r}}} } \\
&= \sum\limits_{\scriptstyle i + j = m \hfill \atop
  \scriptstyle i,j \ge 0 \hfill} {{{\left( { - 1} \right)}^i}\zeta _r^ \star \left( {{{\left\{ 1 \right\}}_i}} \right)\left\{ {{\zeta _{n + r}}\left( {{{\left\{ 1 \right\}}_{j + 1}}} \right) - {\zeta _{n + r}}\left( {{{\left\{ 1 \right\}}_{j + 1}}} \right)} \right\}} \\
  & = \sum\limits_{\scriptstyle i + j = m + 1 \hfill \atop
  \scriptstyle i,j \ge 0 \hfill} {{{\left( { - 1} \right)}^i}\zeta _r^ \star \left( {{{\left\{ 1 \right\}}_i}} \right){\zeta _{n + r}}\left( {{{\left\{ 1 \right\}}_j}} \right)}  - \sum\limits_{\scriptstyle i + j = m + 1 \hfill \atop
  \scriptstyle i,j \ge 0 \hfill} {{{\left( { - 1} \right)}^i}\zeta _r^ \star \left( {{{\left\{ 1 \right\}}_i}} \right){\zeta _r}\left( {{{\left\{ 1 \right\}}_j}} \right)} .\tag{2.10}
\end{align*}
On the other hand, from Lemma 2.5, setting ${a_k} = \frac{1}{k}$ and $n=r$ we get
\[\sum\limits_{\scriptstyle i + j = m \hfill \atop
  \scriptstyle i,j \ge 0 \hfill} {{{\left( { - 1} \right)}^i}\zeta _r^ \star \left( {{{\left\{ 1 \right\}}_i}} \right){\zeta _r}\left( {{{\left\{ 1 \right\}}_j}} \right)}  = 0,\;m \ge 1.\tag{2.11}\]
Hence, combining (2.10) and (2.11) we can prove that the formula (2.8) holds.\hfill$\square$\\
Similarly, by a similar argument as in the proof of Theorem 2.6 with the help of formula (5.2) in the reference \cite{Xu2017}, we obtain the more general theorem.
\begin{thm} For integers $r\geq 0, m,n>1$ and real $p>0$, then
\begin{align*}
&\sum\limits_{1 \le {k_m} <  \cdots  < {k_1} \le n} {\frac{1}{{{{\left( {{k_1} + r} \right)}^p} \cdots {{\left( {{k_m} + r} \right)}^p}}}}  = \sum\limits_{\scriptstyle i + j = m \hfill \atop
  \scriptstyle i,j \ge 0 \hfill} {{{\left( { - 1} \right)}^i}\zeta _r^ \star \left( {{{\left\{ p \right\}}_i}} \right){\zeta _{n + r}}\left( {{{\left\{ p \right\}}_j}} \right)} ,\tag{2.12}\\
&\sum\limits_{1 \le {k_m} \le  \cdots  \le {k_1} \le n} {\frac{1}{{{{\left( {{k_1} + r} \right)}^p} \cdots {{\left( {{k_m} + r} \right)}^p}}}}  = \sum\limits_{\scriptstyle i + j = m \hfill \atop
  \scriptstyle i,j \ge 0 \hfill} {{{\left( { - 1} \right)}^i}{\zeta _r}\left( {{{\left\{ p \right\}}_i}} \right)\zeta _{n + r}^ \star \left( {{{\left\{ p \right\}}_j}} \right)} .\tag{2.13}
\end{align*}
\end{thm}
\begin{re} In fact, in the same way as above, the results of Theorem 2.6 and 2.7 can be extended to the following generalized conclusion.
\begin{align*}
&\sum\limits_{1 \le {k_m} <  \cdots  < {k_1} \le n} {{a_{{k_1} + r}} \cdots {a_{{k_m} + r}}}  = \sum\limits_{\scriptstyle i + j = m \hfill \atop
  \scriptstyle i,j \ge 0 \hfill} {{{\left( { - 1} \right)}^i}{A_i}\left( r \right){{\bar A}_j}\left( {n + r} \right)} ,\\
&\sum\limits_{1 \le {k_m} \le  \cdots  \le {k_1} \le n} {{a_{{k_1} + r}} \cdots {a_{{k_m} + r}}}  = \sum\limits_{\scriptstyle i + j = m \hfill \atop
  \scriptstyle i,j \ge 0 \hfill} {{{\left( { - 1} \right)}^i}{{\bar A}_i}\left( r \right){A_j}\left( {n + r} \right)},
\end{align*}
where ${{A}_m}\left( n \right)$ and ${{\bar A}_m}\left( n \right)$ are defined in Lemma 2.5. It is clear that Lemma 2.5 and Theorem 2.7 are immediate corollaries of Remark 2.1.
\end{re}
\begin{thm} For integers $r\geq 0$ and $m,n\geq 1$, then the following identity holds:
\[h_n^{\left( {r + 1} \right)}\left( m \right) = \left( {\begin{array}{*{20}{c}}
   {n + r}  \\
   r  \\
\end{array}} \right)\sum\limits_{\scriptstyle i + j = m \hfill \atop
  \scriptstyle i,j \ge 0 \hfill} {{{\left( { - 1} \right)}^i}\zeta _r^ \star \left( {{{\left\{ 1 \right\}}_i}} \right){\zeta _{n + r}}\left( {{{\left\{ 1 \right\}}_j}} \right)} .\tag{2.14}\]
\end{thm}
\pf In Lemma 2.3, taking ${x_j} = \frac{1}{{j + r}}$, then we have
\begin{align*}
&{X_n}\left( m \right) = \sum\limits_{j = 1}^n {{{\left( {\frac{1}{{j + r}}} \right)}^m}}  = H_{n + r}^{\left( m \right)} - H_r^{\left( m \right)},\\
&{\zeta _n}\left( {{{\left\{ 1 \right\}}_m}\left| {r + 1} \right.} \right) = \frac{{{{\left( { - 1} \right)}^{m - 1}}}}{m}\sum\limits_{i = 0}^{m - 1} {{{\left( { - 1} \right)}^i}{\zeta _n}\left( {{{\left\{ 1 \right\}}_i}\left| {r + 1} \right.} \right)\left( {H_{n + r}^{\left( {m - i} \right)} - H_r^{\left( {m - i} \right)}} \right)} .\tag{2.15}
\end{align*}
From Lemma 2.1 and formula (2.15) we deduce that
\[h_n^{\left( {r + 1} \right)}\left( m \right) = \left( {\begin{array}{*{20}{c}}
   {n + r}  \\
   r  \\
\end{array}} \right){\zeta _n}\left( {{{\left\{ 1 \right\}}_m}\left| {r + 1} \right.} \right).\tag{2.16}\]
Substituting (2.8) into (2.16) we may easily obtain the desired result. This completes
the proof of Theorem 2.7. \hfill$\square$
\section{Proof of Theorem 1.1}
By replacing $x$ by $n$ and $n$ by $r$ in (1.9), we deduce that
\[\left( {\begin{array}{*{20}{c}}
   {n + r}  \\
   r  \\
\end{array}} \right) = \frac{1}{{r!}}\sum\limits_{l = 1}^{r + 1} {\left[ {\begin{array}{*{20}{c}}
   {r + 1}  \\
   l  \\
\end{array}} \right]{n^{l - 1}}} .\tag{3.1}\]
Therefore, from (2.12) and (3.1) we obtain
\[h_n^{\left( {r + 1} \right)}\left( k \right) = \frac{1}{{r!}}\sum\limits_{l = 1}^{r + 1} {\left[ {\begin{array}{*{20}{c}}
   {r + 1}  \\
   l  \\
\end{array}} \right]{n^{l - 1}}} \sum\limits_{\scriptstyle i + j = k \hfill \atop
  \scriptstyle i,j \ge 0 \hfill} {{{\left( { - 1} \right)}^i}\zeta _r^ \star \left( {{{\left\{ 1 \right\}}_i}} \right){\zeta _{n + r}}\left( {{{\left\{ 1 \right\}}_j}} \right)} .\tag{3.2}\]
Thus, by the definition of $S\left( {k,m;p} \right)$ and (3.2) we can prove (1.7).\hfill$\square$
\section{Proof of Theorem 1.2}
By the definition of multiple harmonic number (1.1), we can find that
\begin{align*}
&{U_{m,r}}\left( p \right) = \sum\limits_{n = 1}^\infty  {\frac{{{\zeta _{n + r}}\left( {{{\left\{ 1 \right\}}_m}} \right)}}{{{n^p}}}}  = \sum\limits_{n = 1}^\infty  {\frac{1}{{{n^p}}}\sum\limits_{k = 1}^{n + r} {\frac{{{\zeta _{k - 1}}\left( {{{\left\{ 1 \right\}}_{m - 1}}} \right)}}{k}} } \\
& = \sum\limits_{n = 1}^\infty  {\frac{{{\zeta _n}\left( {{{\left\{ 1 \right\}}_m}} \right)}}{{{n^p}}}}  + \sum\limits_{k = 1}^r {\sum\limits_{n = 1}^\infty  {\frac{{{\zeta _{n + k - 1}}\left( {{{\left\{ 1 \right\}}_{m - 1}}} \right)}}{{{n^p}\left( {n + k} \right)}}} } \\
&= \zeta \left( {p,{{\left\{ 1 \right\}}_m}} \right) + \zeta \left( {p + 1,{{\left\{ 1 \right\}}_{m - 1}}} \right) + \sum\limits_{{i_1} = 1}^r {\sum\limits_{n = 1}^\infty  {\frac{{{\zeta _n}\left( {{{\left\{ 1 \right\}}_{m - 1}}} \right)}}{{{n^p}\left( {n + {i_1}} \right)}}} }  + \sum\limits_{{i_1} = 1}^r {\sum\limits_{{i_2} = 1}^{{i_1} - 1} {\sum\limits_{n = 1}^\infty  {\frac{{{\zeta _{n + {i_2} - 1}}\left( {{{\left\{ 1 \right\}}_{m - 2}}} \right)}}{{{n^p}\left( {n + {i_1}} \right)\left( {n + {i_2}} \right)}}} } } \\
& =  \cdots \\
& = \zeta \left( {p,{{\left\{ 1 \right\}}_m}} \right) + \zeta \left( {p + 1,{{\left\{ 1 \right\}}_{m - 1}}} \right) + \sum\limits_{j = 1}^m {\sum\limits_{1 \le {i_j} <  \cdots  < {i_1} \le r} {\sum\limits_{n = 1}^\infty  {\frac{{{\zeta _n}\left( {{{\left\{ 1 \right\}}_{m - j}}} \right)}}{{{n^p}\left( {n + {i_1}} \right) \cdots \left( {n + {i_j}} \right)}}} } } .\tag{4.1}
\end{align*}
On the other hand,
we consider the expansion
\[\frac{1}{{\prod\limits_{i = 1}^k {\left( {n + {a_i}} \right)} }} = \sum\limits_{j = 1}^k {\frac{{{A_j}}}{{n + {a_j}}}}\ \ (k\in \N_0; a_i \in \mathbb{C}\setminus\Z^-)\tag{4.2} \]
where
\[{A_j} = \mathop {\lim }\limits_{n \to  - {a_j}} \frac{{n + {a_j}}}{{\prod\limits_{i = 1}^k {\left( {n + {a_i}} \right)} }} = \prod\limits_{i = 1,i \ne j}^k {{{\left( {{a_i} - {a_j}} \right)}^{ - 1}}}.\tag{4.3}\]
Therefore, the equation (4.1) can be written as
\begin{align*}
{U_{m,r}}\left( p \right) =& \zeta \left( {p,{{\left\{ 1 \right\}}_m}} \right) + \zeta \left( {p + 1,{{\left\{ 1 \right\}}_{m - 1}}} \right)\\
& + \sum\limits_{j = 1}^m {\left( {\sum\limits_{l = 1}^j {\prod\limits_{a = 1,a \ne l}^j {{{\left( {{i_a} - {i_l}} \right)}^{ - 1}}} } } \right)\sum\limits_{1 \le {i_j} <  \cdots  < {i_1} \le r} {\sum\limits_{n = 1}^\infty  {\frac{{{\zeta _n}\left( {{{\left\{ 1 \right\}}_{m - j}}} \right)}}{{{n^p}\left( {n + {i_l}} \right)}}} } } .\tag{4.4}
\end{align*}
For $r>0$, we have the partial fraction decomposition
\[\frac{1}{{{n^p}\left( {n + r} \right)}} = \sum\limits_{b = 1}^{p - 1} {\frac{{{{\left( { - 1} \right)}^{b - 1}}}}{{{r^b}}}\cdot\frac{1}{{{n^{p + 1 - b}}}}}  + \frac{{{{\left( { - 1} \right)}^{p - 1}}}}{{{r^{p - 1}}}}\cdot\frac{1}{{n\left( {n + r} \right)}},\tag{4.5}\]
Moreover, from identities (2.1), (2.4) and (2.6), we deduce the following result
\[\sum\limits_{n = 1}^\infty  {\frac{{{\zeta _n}\left( {{{\left\{ 1 \right\}}_{p - 1}}} \right)}}{{n\left( {n + k} \right)}}}  = \frac{1}{k}\left\{ {\zeta \left( p \right) + \zeta _k^ \star \left( {{{\left\{ 1 \right\}}_p}} \right) - \frac{{\zeta _k^ \star \left( {{{\left\{ 1 \right\}}_{p - 1}}} \right)}}{k}} \right\},\;k\in \N,p \in \N/\{1\}.\tag{4.6} \]
Hence, combining (4.4), (4.5) and (4.6), by a simple calculation, we obtain the the desired result. This completes
the proof of Theorem 1.2. \hfill$\square$

Similarly, applying the same arguments as in the proof of formula (4.1), we also deduce a similar result
\begin{align*}
{V_{m,r}}\left( p \right) :&= \sum\limits_{n = 1}^\infty  {\frac{{\zeta _{n + r}^ \star \left( {{{\left\{ 1 \right\}}_m}} \right)}}{{{n^p}}}} \\
& = {\zeta ^ \star }\left( {p,{{\left\{ 1 \right\}}_m}} \right) + \sum\limits_{j = 1}^m {\sum\limits_{1 \le {i_j} \le  \cdots  \le {i_1} \le r} {\sum\limits_{n = 1}^\infty  {\frac{{\zeta _n^ \star \left( {{{\left\{ 1 \right\}}_{m - j}}} \right)}}{{{n^p}\left( {n + {i_1}} \right) \cdots \left( {n + {i_j}} \right)}}} } } .\tag{4.7}
\end{align*}
In fact, we hope to obtain a similar result of Theorem 1.2, but, so far, we have been unable to continue any progress with this sums.
\section{Conclusion}
From \cite{KO2010,OZ2001}, we know that the Aomoto-Drinfel¡¯d-Zagier formula reads
\[\sum\limits_{n,m = 1}^\infty  {\zeta \left( {m + 1,{{\left\{ 1 \right\}}_{n - 1}}} \right){x^m}{y^n} = 1 - \exp \left( {\sum\limits_{n = 2}^\infty  {\zeta \left( n \right)\frac{{{x^n} + {y^n} - {{\left( {x + y} \right)}^n}}}{n}} } \right)} ,\tag{5.1}\]
which implies that for any $m,\ n\in \N$, the multiple zeta value ${\zeta \left( {m + 1,{{\left\{ 1 \right\}}_{n - 1}}} \right)}$ can be
represented as a polynomial of zeta values with rational coefficients, and we have the duality formula
\[\zeta \left( {n + 1,{{\left\{ 1 \right\}}_{m - 1}}} \right) = \zeta \left( {m + 1,{{\left\{ 1 \right\}}_{n - 1}}} \right).\]
In particular, one can find explicit formulas for small weights.
\[\begin{array}{l}
 \zeta \left( {2,{{\left\{ 1 \right\}}_m}} \right) = \zeta \left( {m + 2} \right), \\
 \zeta \left( {3,{{\left\{ 1 \right\}}_m}} \right) = \frac{{m + 2}}{2}\zeta \left( {m + 3} \right) - \frac{1}{2}\sum\limits_{k = 1}^m {\zeta \left( {k + 1} \right)\zeta \left( {m + 2 - k} \right)} . \\
 \end{array}\]
Hence, from formulas (1.11) and (5.1), we see that the sums ${U_{m,r}}\left( p \right)$ can be expressed in terms of series of Riemann zeta function and harmonic numbers. Thus, we show that the Euler-type sums with hyperharmonic numbers
$S\left( {k,r + 1;p} \right)$ can be expressed in terms of zeta values and Stirling numbers of the first kind, for integers $k\in \mathbb{N}$ and $r\in \mathbb{N}_0$ with $r+2\leq p\in \mathbb{N}$.

To conclude, note that it would be useful to be able to extend the approach described above to include other similar and related sums. In particular, it would be very interesting to consider sums of the form
\[{S^ \star }\left( {k,m;p} \right): = \sum\limits_{n = 1}^\infty  {\frac{{H_n^{\left( m \right)}\left( k \right)}}{{{n^p}}}} \;\;\left( {p \ge m + 1} \right),\]
where $H_n^{\left( m \right)}\left( k \right)$ is called the generalized hyperharmonic star numbers, defined by
\[H_n^{\left( m \right)}\left( k \right): = \sum\limits_{\scriptstyle 1 \le {n_{m + k - 1}} \le  \cdots  \le {n_m} \hfill \atop
  \scriptstyle  \le {n_{m - 1}} \le  \cdots  \le {n_1} \le n \hfill} {\frac{1}{{{n_m}{n_{m + 1}} \cdots {n_{m + k - 1}}}}} .\]
Obviously, by the definitions of $H_n^{\left( m \right)}\left( k \right)$ and $h_n^{\left( m \right)}(k)$, we have
\[H_n^{\left( m \right)}\left( 1 \right)=h_n^{\left( m \right)}(1) = \zeta _n^ \star \left( {{{\left\{ 0 \right\}}_{m - 1}},1} \right).\]
Hence $S\left( {1,m;p} \right)={S^ \star }\left( {1,m;p} \right)$.
However, we have been unable, so far, to make any progress with this sums. Unfortunately, it appears that, even in the case $k=2$, the method used in this work gives rise to several complex and intractable summations.

{\bf Acknowledgments.} The authors would like to thank the anonymous
referee for his/her helpful comments, which improve the presentation
of the paper.

 \end{document}